\documentclass[12pt]{article}
\newcommand{\field}[1]{\mathbb{#1}}
\usepackage{graphics,amsmath,amssymb}
\usepackage{natbib}
\bibpunct{[}{]}{;}{a}{,}{,}

\title{Justifying Definitions in Mathematics---Going Beyond Lakatos}
\author{Charlotte Werndl,  c.s.werndl@lse.ac.uk\\ 
Department of Philosophy, Logic and Scientific Method \\London School of Economics}
\date{\normalsize This is a pre-copyedited, author-produced PDF of an article accepted for publication in Philosophia Mathematica following peer review. The definitive publisher-authenticated version ``C. Werndl (2009), Justifying Definitions in Mathematics---Going Beyond Lakatos, Philosophia Mathematica 17, 313-340'' is available online at: http://philmat.oxfordjournals.org/content/17/3/313.abstract}
\begin{document}
\newtheorem{definition}{Definition}
\newtheorem{lemma}{Lemma}
\newtheorem{theorem}{Theorem}
\newtheorem{remark}{Remark}
\maketitle
\begin{abstract}
This paper addresses the actual practice of justifying definitions in mathematics. First, I introduce the main account of this issue, namely Lakatos's proof-generated definitions. Based on a case study of definitions of randomness in ergodic theory, I identify three other common ways of justifying definitions: natural-world-justification, condition-justification and redundancy-justification. Also, I clarify the interrelationships between the different kinds of justification. Finally, I point out how Lakatos's ideas are limited: they fail to show that various kinds of justification can be found and can be reasonable, and they fail to acknowledge the interplay between the different kinds of justification.
\end{abstract}
\tableofcontents
\newpage
\section{Introduction}

Mathematical practice suggests that mathematical definitions are not arbitrary: for definitions to be worth studying there have to be good reasons. Moreover, definitions are often regarded as important mathematical knowledge [cf.~\citeauthor{Tappenden2008a}, \citeyear{Tappenden2008a} and \citeyear{Tappenden2008b}]. Reasoning and knowledge are classical philosophical issues; hence reflecting on the reasons given for definitions is philosophically relevant.

These considerations motivate the guiding question of this paper:
\textit{in what ways are definitions in mathematics justified, and are these kinds of justification reasonable?} By a justification of a definition we mean a reason provided for the definition.
We will concentrate on \textit{explicit} definitions, which introduce a new expression by stipulating that it be semantically equivalent to the definiens consisting of already-known expressions. We won't deal with their complement, \textit{implicit} definitions, which assign meaning to expressions by imposing constraints on how to use sentences (or other longer expressions) containing them \cite[p.~97]{Brown1999}.

Generally, attempting to justify definitions is reasonable: as we will see, if definitions were not justified, the mathematics involving these definitions would be much less meaningful to us than mathematics involving definitions which were justified. Thus given our limited resources, it is better to concentrate on definitions which we can justify.\footnote{
What this means for the ontology of mathematical definitions depends on the ontology adopted: platonists may hold that both definitions which we can and cannot justify are real, but it is better to concentrate on the former for pragmatic reasons; constructivists may hold that only definitions which we can justify are real.}

When a mathematician formulates a definition she or he has not known before, we speak of a formulation of the definition. The way a formulation of a definition is guided usually corresponds to the way the definition is justified when it is formulated. Thus all that will be said about the justification of definitions has a natural counterpart in terms of the guidance of the formulation of definitions. Since the guidance of the formulation of definitions derives from the justification, the latter is the main issue, and in what follows we will focus on the justification of definitions.\footnote{Strictly speaking, the justification and the guidance of formulation are conceptually distinct. For instance, it could be that a definition which captures an important preformal idea was randomly formulated by a computer; then there was no way the formulation of the definition was guided, but there is a convincing initial justification.}

In this paper, in section \ref{SOTA}, we will first discuss the state of the art of theorising about the actual mathematical practice of how definitions are justified in articles and books. There is hardly any philosophical discussion on this issue apart from Lakatos's ideas on proof-generated definitions, and hence we will concentrate on them. While Lakatos's ideas are important, this paper aims to show how they are limited. My criticism of Lakatos starts from a case study of notions of randomness in ergodic theory, which will be introduced in section \ref{CS}. Based on this case study, in section \ref{chaos} I will introduce three other  ways in which definitions are commonly justified: natural-world justification, condition-justification and redundancy-justification; two of them, to my knowledge, have not been discussed before. In section \ref{interr} I will clarify the interrelationships between the different kinds of justification, an issue which also has not been addressed before. In particular, I argue that in different arguments the same definition can be justified in different ways. Finally, in section \ref{Lakatos} I point out how Lakatos's ideas are limited: his ideas fail to show that often, as for notions of randomness in ergodic theory, various kinds of justification are found and that various kinds of justification can be reasonable. Furthermore, they fail to acknowledge the interplay between the different kinds of justification.

This research is in the spirit of `phenomenological philosophy of mathematics' as recently characterised by \citeauthor{Larvor2001}~[\citeyear{Larvor2001}, pp.~214--215] and \citeauthor{Leng2002}~[\citeyear{Leng2002}, pp.~3--5]: it looks at mathematics `from the inside' and on this basis asks philosophical questions.

\section{Lakatos's Proof-Generated Definitions}\label{SOTA}
In the relatively recent literature \citeauthor{Larvor2001}~[\citeyear{Larvor2001}, p.~218] at least mentions the importance of researching the justification of mathematical definitions. \citeauthor{Corfield2003}~[\citeyear{Corfield2003}, chapter~9]  discusses the related issue of what makes fundamental concepts but does not provide conceptual reflection on our question. \citeauthor{Tappenden2008a}~[\citeyear{Tappenden2008a} and \citeyear{Tappenden2008b}] treats the related issues of naturalness of definitions and how to decide between different definitions. In our context \citeauthor{Tappenden2008a}'s [\citeyear{Tappenden2008a}] conclusion is relevant: namely, that judgments about definitions mainly depend not on the rules of logic but on detailed knowledge about the mathematics involved. Furthermore, several philosophers have argued that mathematical definitions should capture a valuable preformal idea [cf.~\citeauthor{Brown1999}, \citeyear{Brown1999}, p.~109].

Apart from this, the main philosopher who has written on our guiding question in the light of mathematical practice is \citeauthor{Lakatos1976}~[\citeyear{Lakatos1976}]. Lakatos develops an approach of informal mathematics, which includes an account of mathematical progress called \textit{proofs and refutations}. Most importantly, Lakatos is also concerned with how definitions are justified. His key idea is the notion of a \textit{proof-generated definition}. Here his main example are definitions of polyhedron which are justified because they are needed to make the proof of the Eulerian conjecture work: viz.\ that for every polyhedron the number of vertices minus the number of edges plus the number of faces equals two ($V-E+F=2$).

What is a proof-generated definition? Unfortunately, Lakatos does not state exactly what he means by this. Clearly, mathematical definitions justified in any way are eventually involved in proofs. Therefore, the trivial idea that definitions are justified because they are involved in proofs cannot be what interested Lakatos.

To find out more, consider the \citeauthor{Caratheodory1914} definition of measurable sets, another proof-generated definition Lakatos discusses. The mathematician \citet[p.~44]{Halmos1950} remarks on this definition: ``The greatest justification of this apparently complicated concept is, however, its possibly surprising but absolute complete success as a tool of proving the extension theorem''. \citeauthor{Lakatos1976}~[\citeyear{Lakatos1976}, p.~153] comments:
\begin{quote}
as we learn from the second part [Halmos's remark above], this concept is a proof-generated concept in \citeauthor{Caratheodory1914}'s theorem about the extension of measures [...]. So whether it is intuitive or not is not at all interesting: its rationale lies not in its intuitiveness but in its proof-ancestor.
\end{quote}
This quote and the rest of the discussion of proof-generated definitions suggests that a \textit{proof-generated definition is a definition which is needed in order to prove a specific conjecture regarded as valuable} \cite[pp.~88--92, pp.~127--133, pp.~144--54]{Lakatos1976}. This idea is also hinted at by \citeauthor{Polya1954}~[\citeyear{Polya1954}, p.~148]. The final theorems which involve proof-generated definitions often, but not always, result  from a series of trials and revisions.

\citeauthor{Lakatos1976}~[\citeyear{Lakatos1976}, pp.~33--50, p.~127] rightly argues that \textit{lemma-incorporation} produces proof-generated definitions: assume that a conjecture, known not to hold for all objects of a domain, should be established. Then if conditions which are needed in order to prove the conjecture are identified, i.e.~lemmas are incorporated, proof-generated definitions arise. For instance, consider the conjecture that the limit function of a convergent sequence of continuous functions is continuous. This conjecture can be proven if `convergent' is understood as uniformly convergent but not if it is understood as the more obvious, weaker pointwise convergent; hence the definition of uniformly convergent is proof-generated \cite[pp.~144--146]{Lakatos1976}.

\citeauthor{Lakatos1976}~[\citeyear{Lakatos1976}, pp.~90--92, p.~128, pp.~148--149, p.~153] thinks that for his examples of proof-generated definitions the justification was \textit{reasonable} because the corresponding conjectures are valuable. Generally, if the conjecture is mathematically valuable, proof-generation is a reasonable kind of justification.\footnote{For the proof-generated definitions discussed in \cite{Lakatos1976} and in this it is argued why the conjectures are valuable. Yet answering the question of what constitutes valuable conjectures \textit{at a general level} would require further research.} A proof-generated definition can be regarded as providing knowledge since it answers the question of which notion is needed to prove a specific conjecture.

\citeauthor{Lakatos1976}~[\citeyear{Lakatos1976}, pp.~14--33, pp.~83--87] also discusses four other ways of justifying definitions. Imagine that counterexamples are presented to a conjecture of interest, and that
the conjecture is defended by claiming that these are no ``real'' counterexamples because a definition in the conjecture has been wrongly understood. Properly understood, it is argued, the definition excludes a class of objects which includes the alleged counterexamples, where the exclusions are made independent of any proof of the conjecture (and thus it is unknown whether the conjecture indeed holds true for the definition). Then the definition is justified via \textit{monster-barring}.
The second kind of justification is \textit{exception-barring}.  Here the definition is defended by excluding, with the extant definition, a class of objects which are counterexamples to the conjecture; again, this is independent of any proof of the conjecture.\footnote{Contrary to exception-barring, in the case of monster-barring
it is denied that the counterexamples are actual counterexamples. This is how monster-barring differs from exception-barring.}
The third kind of justification is \textit{monster-adjustment}. Here the definition is defended by reinterpreting, independent of any proof of the conjecture, the terms of the extant definition such that counterexamples to the conjecture are no longer counterexamples any longer.
The fourth and final kind of justification is \textit{monster-including}.  Here the definition is defended by extending the definition to include a new class of objects; this class of objects is defined using properties which are shared by examples for
which the conjecture holds true; and again, this is independent of any
proof of the conjecture.

Monster-barring, exception-barring and monster-adjustment are all ways of dealing with counterexamples to conjectures. And I agree with Lakatos that for this purpose they are inferior to proof-generation because they do not take into account how the conjectures are proved; and therefore, it is even unclear whether the conjecture is true for the definition under consideration. Monster-including is a way of generalising conjectures. Yet again, since it neglects how conjectures are proved, I agree with Lakatos that for this purpose it is inferior to proof-generation. Furthermore, Lakatos thought that any of these kinds of justification were applied only because the better way of justifying definitions, namely with proof-generation, was not known [\citeauthor{Lakatos1976}~\citeyear{Lakatos1976}, pp.~14--42, pp.~136--140]. Because of their inadequacies and since they play no role in our case study, I shan't say any more about these kinds of justification in this paper.

Unfortunately, \citeauthor{Lakatos1976}~[\citeyear{Lakatos1976}] never explicitly states how widely he thinks that his ideas on proof-generated definitions apply. General claims such as
\begin{quote}
Progress indeed replaces \textit{naive classification} by [...]\,proof-generated [...]\linebreak[3]classification. [...]\,\textit{Naive conjectures and naive concepts are superseded by improved conjectures (theorems) and concepts (proof-generated} [...]\,\textit{concepts) growing out of the method of proofs and refutations} [\citeauthor{Lakatos1976},~\citeyear{Lakatos1976}, pp.~91--92; see also p.~144].
\end{quote}
suggest that mathematical definitions should be, and after the discovery of proof-generation are generally proof-generated, and some have interpreted him as saying this \cite[pp.~110--111]{Brown1999}.
However, as \citeauthor{Larvor1998}~[\citeyear{Larvor1998}] has pointed out, Lakatos stresses in his PhD-thesis, on which his \citeyear{Lakatos1976} book is based, that his account of informal mathematics does not apply to all of mathematics.
What is clear is that Lakatos thought that there are many mathematical subjects with some proof-generated definitions and that there are many mathematical subjects with some definitions which should be proof-generated.\footnote{Of course, the question remains what a `mathematical subject' is; I will say more about this later [cf.~section~\ref{subject}].} Maybe Lakatos also believed something stronger, and this would explain his strong claims such as in the above quote, namely that there are many subjects where proof-generation should be the \textit{sole} important way in which definitions are justified; and that there are many subjects created after the discovery of proof-generation where proof-generation is the \textit{sole} important way in which definitions are justified. In what follows, I will show in which ways Lakatos's ideas on justifying definitions are limited; and for this it won't matter much whether or not he endorsed the stronger claim.

\citeauthor{Corfield1997}~[\citeyear{Corfield1997}, pp.~111--115] argues that Lakatos did not think that his account of informal mathematics, which includes his ideas on justifying definitions, extends to established branches of mathematics of ``the twentieth century and up to the present day''. Yet \citeauthor{Corfield1997}'s claim is implausible. \citeauthor{Lakatos1976}~[\citeyear{Lakatos1976}, p.~5, pp.~152--154] states that his ideas on informal mathematics apply to modern metamathematics and to \citeauthor{Caratheodory1914}'s~[\citeyear{Caratheodory1914}] investigations on measurable sets. And substantial parts of established mathematics of the twentieth century up to the present day are not more formalised than that mathematics: e.g.~ergodic theory, which will be relevant later. Thus Lakatos indeed thought that his ideas could apply to substantial parts of established branches of mathematics of the twentieth century up to the present day. But I agree with \citeauthor{Corfield1997}'s [\citeyear{Corfield1997}] main point that Lakatos failed to see that his ideas are also relevant for highly formalised mathematics. For this reason this paper is not restricted to informal mathematics.

This discussion highlights that there is little work on the actual practice of how definitions are justified in articles and books. Furthermore, although Lakatos's account of proofs and refutations has been challenged \citep{Corfield1997,Leng2002}, his ideas on proof-generated definitions have hardly been criticised. My contribution on the guiding question and my criticism of Lakatos's ideas on justifying definitions will be based on a case study of notions of randomness in ergodic theory. Let me now introduce this case study.

\section{Case Study: Randomness in Ergodic Theory}\label{CS}

My case study is on \textit{notions of randomness in ergodic theory}. Ergodic theory originated from work in statistical mechanics, in particular Boltzmann's kinetic theory of gases. Boltzmann's work relied on the assumption that the time-average of a function equals its space average, but no acceptable argument  was provided for this. Generally, the possible random motion of classical systems was a constant theme in statistical mechanics. Ergodic theory arose in the early 1930s when von Neumann and Birkhoff proved the famous mean and pointwise ergodic theorem, respectively. Among other things, they found that ergodicity was the sought-after concept guaranteeing the equality of time and space averages for almost all states of the system. Motivated by these results, an investigation into the random behaviour of classical systems began. Of particular importance here was the study of randomness by a group of mathematicians around Kolmogorov in Russia. From the 1960's onwards, ergodic theory became prominent, and was further developed, as a mathematical framework for studying chaotic behaviour, i.e.~unpredictable and random behaviour of deterministic systems. Overall, ergodic theory had less impact on statistical mechanics than expected, partly because of the doubts, and the difficulty of proving, that the relevant systems are ergodic. But it developed into a discipline with its own internal problems and had, and continues to have, considerable impact on probability theory and chaos research [\citeauthor{AubinDahan2002}, \citeyear{AubinDahan2002}; \citeauthor{Mackey1974}, \citeyear{Mackey1974}].

Why do notions of randomness in ergodic theory constitute a valuable case study? First, several of Lakatos's assertions, e.g.~that mathematics is driven by counterexamples, have been criticised in the following way: while they may be correct for older mathematics, they do not hold true for twentieth century mathematics \cite[p.~10]{Leng2002}. As also \citeauthor{Lakatos1976}~[\citeyear{Lakatos1976}, pp.~136--140] suggests, how definitions are justified may depend on when they were formulated because reasoning changes with the advancement of mathematics. To ensure that claims on the justification of definitions escape the criticism of not applying to twentieth century mathematics, I choose mathematics, like ergodic theory, which was created in the twentieth century. Second, concerning the justification of definitions the picture for notions of randomness in ergodic theory appeared different to that proposed by Lakatos, and this picture seemed prevalent in mathematics.

As widely acknowledged, the main notions of randomness in ergodic theory are [cf.~\citeauthor{Berkovitzetal2006}, \citeyear{Berkovitzetal2006}; \citeauthor{Sinai2000}, \citeyear{Sinai2000}, p.~21, pp.~41--46; \citeauthor{Walters1982}, \citeyear{Walters1982}, pp.~39--41, pp.~86--87, pp.~105--107]:\label{list}\\\\
\noindent \textit{weak mixing (three versions),} strong mixing (two versions), Kolmogorov-mixing, Kolmogorov-system, \textit{Bernoulli-system (two versions)},  \textit{Kolmogorov-Sinai entropy}.\\

\noindent I studied how the definitions listed above are justified and whether
they are reasonably justified. I also examined the way the definitions were initially justified.\footnote{I did not investigate the use of these definitions elsewhere in mathematics. The main reason for such an investigation would have been to understand how the justification of definitions varies in different contexts. Yet I could find out about this by considering only how definitions were initially justified and later justified as notions of randomness. Going further would have required an enormous amount of work without considerable gain.}

In the remaining sections of this paper the insights on the justification of definitions which derive from this case study will be presented. The definitions of the above list which are italicised will be discussed in detail. A detailed investigation of them will suffice to illustrate these insights. Hence for the remaining listed definitions I will just state how they are justified. Basic knowledge of measure theory will suffice to understand the mathematics that follows. Yet the claims can also be grasped without understanding the relevant mathematics if close attention is paid to the verbal commentary on the definitions.

First of all, the unit of analysis in ergodic theory, namely measure-theoretic dynamical systems, needs to be introduced.
Generally, dynamical systems are mathematical models consisting of a phase space, the set of all possible states of the system, and an evolution equation which describes how solutions evolve in phase space. Dynamical systems usually model natural systems.

There are discrete and continuous dynamical systems. For discrete systems time increases in discrete steps. Continuous systems involve a continuous time parameter; they typically arise from differential equations.
Only for one definition of randomness in ergodic theory there is a considerable difference between the discrete and continuous version, and here I will discuss both versions. All other notions of randomness are essentially the same for discrete and continuous time, and hence, for simplicity, we confine our attention to the discrete notions.

Ergodic theory is concerned with \textit{measure-theoretic dynamical systems}. A discrete measure-theoretic dynamical system is a quadruple $(X,\Sigma,\mu,T)$ where $X$ is a set (phase space), $\Sigma$ is a $\sigma$-algebra on $X$, $\mu$ is a measure with $\mu(X)=1$ and $T:X\rightarrow X$ (evolution equation) is bijective measurable map which is measure-preserving, i.e.~$\mu(T(A))=\mu(A)$ for all $A\in\Sigma$. A continuous measure-theoretic dynamical system is a quadruple $(X,\Sigma,\mu,T_{t})$ where
$X$ is a set (phase space), $\Sigma$ is a $\sigma$-algebra on $X$, $\mu$ is a measure with $\mu(X)=1$, and $T_{t}:X\rightarrow X$, $t\in\field{R}$, (evolution equation) is a group of measurable and measure-preserving bijective maps such that additionally the function $\Phi(x,t):=T_{t}(x)$ is measurable on $X\times \field{R}$  [cf.~\citeauthor{Cornfeldetal1982}~\citeyear{Cornfeldetal1982}, pp.~3--6].

With this background we can now discuss the kinds of justification which occur in my case study. They illustrate that not only proof-generation is important.

\section{Kinds of Justification of Definitions}\label{chaos}
\subsection{Natural-World-Justification}

I claim, first, that definitions in my case study are frequently justified because \textit{they capture a preformal idea regarded as valuable for describing or understanding the natural world}. Here I will speak of \textit{natural-world-justified} definitions. Natural-world-justified definitions are a special case of the general idea discussed in the literature that mathematical definitions should
capture a valuable preformal idea [cf.~\citeauthor{Brown1999}, \citeyear{Brown1999}, p.~109].

If the preformal idea is valuable for describing or understanding the natural world, natural-world-justification is reasonable. It is important to realise that natural-world-justification does not mean that there is a `best' definition of a vague idea. There can be several different definitions expressing a vague idea without a clearly `best' one. Natural-world-justified definitions can be regarded as providing knowledge in the following sense: they are a possible formalisation of a preformal idea which is valuable.

Many definitions of the list of notions of randomness [cf.~section \ref{CS}] are natural-world-justified: we will now discuss weak mixing (one version), Bernoulli-system (one version) and the Kolmogorov-Sinai entropy in detail. For illustrating natural-world-justification, it would suffice to consider the Kolmogorov-Sinai entropy. The discussion of the remaining two definitions is crucial in order to provide the necessary background for the next sections.
Moreover, all versions of strong mixing [\citeauthor{Berkovitzetal2006}, \citeyear{Berkovitzetal2006}, p.~676; \citeauthor{Hopf1932a}, \citeyear{Hopf1932a}, p.~205] and Kolmogorov-mixing \cite[p.~66]{Sinai1963} are  natural-world-justified.

\subsubsection{Weak Mixing}

\begin{definition}\label{wmsA}
$(X,\Sigma,\mu,T)$ is weakly mixing iff for all $A,B\in\Sigma$ there is a $P\subseteq\field{N}$ of density zero such that
\begin{displaymath}
\lim_{n\rightarrow\infty,\,\,n\notin P}\mu(T^{n}(A)\cap B)=\mu(A)\mu(B),
\end{displaymath} where $P\subseteq\field{N}$ is of density zero iff $\lim_{n\rightarrow\infty}\#(P\cap\{i\,|\,i\leq n,\,\,i\in\field{N}\})/n=0$.\end{definition}

For $(X,\Sigma,\mu,T)$ and $A\in\Sigma$ define $A^{t}$ as the  \textit{event} that the system's state is in $A$ at time $t$, and let $p(...)$ denote the probability of events. Because time is discrete, we can denote time points by $t_{n},\,\,n\in\field{Z}$. Assume, as it is often done, that the measure can be interpreted as time-independent probability: $p(A^{t})=\mu(A)$ and $p(A^{t}\&B^{t})=\mu(A\cap B)$ for all $t$ and all $A,B\in\Sigma$ [see \citeauthor{Werndl2009}, 2009, for more on this]. Then it follows that for all $t_{m},t_{n},\,\,t_{m}\geq t_{n}$, and all $A,B\in\Sigma$:\footnote{$T^{m-n}(A)$ is the evolution of $A$ forward in time from $t_{n}$ to $t_{m}$, containing exactly those points that are in $A$ at time $t_{n}$. Consequently, $T^{m-n}(A)\cap B$ consists of exactly those points which pass $A$ at time $t_{n}$ and go through $B$ at time $t_{m}$, i.e.\ for which $A^{t_{n}}\&B^{t_{m}}$ is true. Therefore, $p(A^{t_{n}}\&B^{t_{m}})=\mu(T^{m-n}(A)\cap B)$.}\begin{equation}\label{EventsSystems}
p(A^{t_{n}}\&B^{t_{m}})=\mu(T^{m-n}(A)\cap B).\end{equation}

By writing out the definition of limit, the definiens of \textit{Definition~\ref{wmsA}} says that for any $A,B\in\Sigma$ and any $\varepsilon>0$ there is a $n_{0}\in\field{N}$ and a set $P$ of density zero with $|\mu(T^{n}(A)\cap B)\!-\!\mu(A)\mu(B)|\!<\!\varepsilon$ for all $n\geq n_{0},\,\,n\notin P$. Equation (\ref{EventsSystems}) thus tells us that \textit{Definition~\ref{wmsA}} captures the following idea of randomness: given an arbitrary level of precision $\varepsilon>0$ any event is approximately independent of almost any event that is sufficiently past. Independence is understood here as in probability theory. This randomness might apply, e.g.~to systems in meteorology and make it hard to predict them.

Von Neumann [\citeyear{Neumann1932}, p.~591, p.~594] lists the main statistical properties of classical systems discussed in ergodic theory and that time. In this context he remarks that \textit{Definition~\ref{wmsA}} captures the preformal idea of approximate independence of almost all events explained above. Thus he argues that it is natural-world justified. This justification grew in importance with the rise of chaos research in the 1960s \cite[cf.][p.~688]{Berkovitzetal2006}. It also appears in a few standard books on ergodic theory \cite[e.g.][p.~45]{Walters1982}, although in books often no justification is provided for weak mixing [e.g.~\citeauthor{Arnold1968}, \citeyear{Arnold1968}, pp.~21-22; \citeauthor{Cornfeldetal1982}, \citeyear{Cornfeldetal1982}, pp.~22--23; \citeauthor{Sinai2000}, \citeyear{Sinai2000}, p.~21].

Especially before the rise of chaos research weak mixing appears to be mostly not naturally-world justified. We will see this in section~\ref{condition}, where we will also discuss the key contexts in which weak mixing was introduced.

The next definition relates to the important topic of equivalence of measure-theoretic systems.
\subsubsection{Discrete-time Bernoulli-system}\label{Bern}
The idea of an infinite sequence of independent trials of an $N$-sided die is a very old one.
\cite{Kolmogorov1933} gave the modern measure-theoretic formulation of probability theory and laid the foundations for the theory of stochastic processes \cite[pp.~230--233]{Plato1994}. In this modern framework an independent process, i.e.~a doubly-infinite sequence of independent rolls of an $N$-sided die where the probability of obtaining $k$ is $p_{k}$, $k\!\in\!\bar{N}:=\{1,\ldots,N\}$, with $\sum_{k=1}^{N}p_{k}\!=\!1$, is modeled as follows. Let $X$ be the set of all bi-infinite sequences $(\ldots x_{-1},x_{0},x_{1}\ldots)$ with $x_{i}\in\bar{N}$, corresponding to the possible outcomes of an infinite sequence of independent trials. Let $\Sigma$ be the set of all sets of infinite sequences to which probabilities can be assigned, and let $\mu$ be the probability function on $\Sigma$.\footnote{In detail: $\Sigma$ is the $\sigma$-algebra generated by the cylinder-sets\begin{displaymath}
C^{k_{1}...k_{n}}_{i_{1}...i_{n}}\!=\!\{x\in X\,|\,x_{i_{1}}\!=\!k_{1},...,x_{i_{n}}\!=\!k_{n},\,i_{j}\!\in\!\field{Z},\,i_{1}\!<...<\!i_{n},\,k_{j}\in\bar{N},\,1\!\leq j\!\leq n\}.\end{displaymath}
The sets have probability $\bar{\mu}(C^{k_{1}...k_{n}}_{i_{1}...i_{n}})=p_{k_{1}}p_{k_{2}}\ldots p_{k_{n}}$
since the outcomes are independent. $\mu$ is defined as the unique extension of $\bar{\mu}$ to a measure on $\Sigma$.} The shift\begin{equation}T:X\rightarrow X\,\,\,\,\,\,\,T((\ldots x_{i}\ldots))=(\ldots x_{i+1}\ldots)\end{equation}
is easily seen to be measurable and measure-preserving. $(X,\Sigma,\mu,T)$ is called a \textit{Bernoulli-shift}.

In one of the first papers on ergodic theory \cite{Neumann1932} introduced the fundamental idea of equivalence of measure-theoretic systems. He developed the definition of \textit{isomorphic systems} to capture this idea \cite[p.~833]{Sinai1989},\footnote{$(X_{1},\Sigma_{1},\mu_{1},T_{1})$ is isomorphic to $(X_{2},\Sigma_{2},\mu_{2},T_{2})$ iff there are measurable sets $\bar{X}_{i}\subseteq X_{i}$ with $\mu_{i}(X_{i}\setminus\bar{X}_{i})=0$ and $T_{i}\bar{X}_{i}\subseteq\bar{X}_{i}\,\,(i=1,2$), and there is a bijection $\phi:\bar{X}_{1}\!\rightarrow \!\bar{X}_{2}$ such that
(i) $\phi(A)\!\in\!\Sigma_{2}$ for all $A\!\in\!\Sigma_{1},A\subseteq \bar{X_{1}}$, and $\phi^{-1}(B)\in\Sigma_{1}$ for all $B\in\Sigma_{2},B\subseteq \bar{X_{2}}$; (ii) $\mu_{2}(\phi(A))=\mu_{1}(A)$ for all $A\in\Sigma_{1},\,A\subseteq\bar{X}_{1}$;
(iii) $\phi(T_{1}(x))=T_{2}(\phi(x))$ for all $x\in\bar{X}_{1}$.} and called for a classification of systems up to isomorphism.

Consequently, we see that the following definition captures the idea of systems which are equivalent to a system describing an independent process, e.g.~throwing a die:
\begin{definition}\label{BerSys}
$(X,\Sigma,\mu,T)$ is a Bernoulli-system iff it is isomorphic to a Bernoulli shift.
\end{definition}
In many articles \textit{Definition~\ref{BerSys}} is natural-world-justified in that way [\citeauthor{Ornstein1989}, \citeyear{Ornstein1989}, p.~4; \citeauthor{Rohlin1960}, \citeyear{Rohlin1960}, p.~5]. Walter's [\citeyear{Walters1982}, p.~107; see also \citeauthor{Ornstein1974}, \citeyear{Ornstein1974}, p.~4] comment
\begin{quote}
Since a Bernoulli shift is really an independent identically distributed stochastic process indexed by the integers we can think of a \{Bernoulli-system\} as an abstraction of such a stochastic process.\footnote{Square brackets indicate that the original notation has been replaced by the notion used in this paper.  I will use this convention throughout.}
\end{quote}
shows that this justification is found in standard books on ergodic theory too. Yet some books do not provide any justification for \textit{Definition~\ref{BerSys}}  [e.g.~\citeauthor{Shields1973}, \citeyear{Shields1973}, p.~5].

Clearly, the Bernoulli-shifts given by choices of $N$ and, for each $N$, the choices of $p_{1},\ldots p_{N}$ are Bernoulli-systems. We will say later more about when Bernoulli-shifts are isomorphic.

The next definition illustrates that a definition can be both preformal-justified and proof-generated.

\subsubsection{Kolmogorov-Sinai Entropy}\label{KSentropy}
$\alpha=\{a_{1},\ldots,a_{k}\}$, $k\geq 1$, is a partition of a system $(X,\Sigma,\mu,T)$ iff $a_{i}\in\Sigma$, $a_{i}\cap a_{j}=\emptyset$ for $i\neq j$, $1\leq i,j\leq k$, and $X=\bigcup_{i=1}^{k}a_{i}$. Clearly, $T^{n}\alpha:=\{T^{n}a_{1},\ldots,T^{n}a_{k}\}$, $n\in\field{Z}$, is also a partition.

Dynamical systems and information theory can be connected as follows: each $x\in X$ produces, relative to a partition $\alpha$, an infinite string of symbols $x_{0}x_{1}x_{2}\ldots$ in an alphabet of $k$ letters via the coding $x_{j}=\alpha_{i}$ iff $T^{j}(x)\in\alpha_{i},\,\,j\geq 0$. Interpreting the system $(X,\Sigma,\mu,T)$ as the source, the output of the source are these strings $x_{0}x_{1}x_{2}\ldots$.
If the measure is interpreted as time-independent probability, $H(\alpha,T):=$\begin{equation}\lim_{n\rightarrow\infty}1/n\!\!\!\!\!\!\!\!\!\!\!\!\!\sum_{i_{j}\in\{1,\ldots,k\},0\leq j\leq n-1}\!\!\!\!\!\!\!\!\!\!\!\!\!\!\!\!\!-\mu(\alpha_{i_{0}}\cap T\alpha_{i_{1}}\ldots\cap T^{n-1}\alpha_{i_{n-1}})\log(\mu(\alpha_{i_{0}}\ldots\cap T^{n-1}\alpha_{i_{n-1}}))\end{equation}
measures the average information which the system produces per step relative to $\alpha$ as time goes to infinity \cite[pp.~233--240]{Petersen1983}. Now:
\begin{definition}\label{KSE}
$h(T):=\sup_{\alpha}\{H(\alpha,T)\}$ is the Kolmogorov-Sinai entropy of the system $(X,\Sigma,\mu,T)$.
\end{definition}
It measures the highest average amount of information that the system can produce per step relative to a coding, and a positive entropy indicates that information is produced.

Having worked for several years on information theory, \cite{Kolmogorov1958} was the first to apply information-theoretic ideas to ergodic theory. He introduced a definition of entropy only for what is nowadays called Kolmogorov-systems. Based on Kolmogorov's work, \cite{Sinai1959} introduced a different notion of entropy which applies to all systems, the now canonical \textit{Definition~\ref{KSE}}. Sinai also proved---a big surprise at that time---that automorphisms on the torus have positive entropy and thus are random because they produce information. Kolmogorov and Sinai were motivated by finding a concept which characterises the amount of randomness of a system \citep{Sinai2007, Werndl2009b}.
More specifically, as Halmos [\citeyear{Halmos1961}, p.~76; cf.~\citeauthor{Sinai1959}, \citeyear{Sinai1959}] explains: ``Intuitively speaking, the entropy $h(T)$ is the greatest quantity of information obtainable about the universe per day [i.e.~step] by repeated performances of experiments with a finite [...] number of outcomes''. Hence \textit{Definition~\ref{KSE}} is natural-world-justified by capturing the idea of the average amount of information produced per step explained above.

Also in some standard books on ergodic theory \textit{Definition~\ref{KSE}} is natural-world-justified in this way [\citeauthor{Billingsley1965}, \citeyear{Billingsley1965}, p.~63; \citeauthor{Petersen1983}, \citeyear{Petersen1983}, pp.~233--240]. It should, however, be mentioned that in books \textit{Definition~\ref{KSE}} is often not justified at all [e.g.\ \citeauthor{Arnold1968}, \citeyear{Arnold1968}, pp.~35--50; \citeauthor{Cornfeldetal1982}, \citeyear{Cornfeldetal1982}, pp.~246--257; \citeauthor{Sinai2000}, \citeyear{Sinai2000}, pp.~40--43].

Interestingly, the Kolmogorov-Sinai entropy is \textit{also proof-generated}. And it is the only notion of randomness in ergodic theory [cf.~section~\ref{CS}] which is proof-generated. The central internal problem of ergodic theory is the following: which systems are isomorphic [cf.~subsection~\ref{Bern}]? Using the so-called Koopman formalism, dynamical systems can be investigated from a spectral-theoretic viewpoint. Systems which are equivalent from this viewpoint are said to be \textit{spectrally-isomorphic}. In the 1950s it was known that systems with discrete spectrum are isomorphic if and only if they are spectrally-isomorphic and that this is not so for systems with mixed spectrum. Most importantly, however, is the case of a \textit{continuous spectrum} since dynamical systems typically have this property~\cite[pp.~27--32]{Arnold1968}.\footnote{Systems have continuous spectrum iff their only eigenfunctions are constant functions.} For systems with continuous spectrum, e.g.~Bernoulli-systems, the conjecture emerged that spectrally-isomorphic systems are not always isomorphic, but the problem resisted solution.

\cite{Kolmogorov1958} and \cite{Sinai1959} were motivated by making progress about this conjecture \cite[p.~834--836]{Sinai1989}, and Kolmogorov's [\citeyear{Kolmogorov1958}] main result is that this conjecture is true. As hinted at by \citet[pp.~1--2, p.~8]{Rohlin1960}, the Kolmogorov-Sinai entropy can be justified as being precisely the definition which is needed to prove that conjecture, i.e.~it is proof-generated. The argument, which goes back to Kolmogorov's work, is as follows: isomorphic system have the same Kolmogorov-Sinai entropy.
Now look at Bernoulli-shifts, whose Kolmogorov-Sinai entropy $\sum_{i}p_{i}\log(p_{i})$ takes a continuum of different values.
Since all Bernoulli-shifts are spectrally-isomorphic, there is a continuum of systems being spectrally-isomorphic but not isomorphic.

Billingsley's [\citeyear{Billingsley1965}, p.~65] comment
\begin{quote}
It is essential to understand the difference between $H(\alpha,T)$ and $h(T)$ and why the latter is introduced. If the entropy of $T$ were taken to be $H(\alpha,T)$ for some ``naturally'' selected $\alpha$ [...], then it would be useless for the isomorphism problem.
\end{quote}
shows that the justification of \textit{Definition~\ref{KSE}} as being proof-generated made it into standard books on ergodic theory too [see also \citeauthor{Petersen1983}, \citeyear{Petersen1983}, p.~227, p.~246].

Let us turn to the second kind of justification I have identified.

\subsection{Condition-Justification}\label{condition}
I claim that another kind of justification abounds in my case study: \textit{a definition is justified by the fact that it is equivalent in an allegedly natural way to a previously specified condition which is regarded as mathematically valuable}. We speak here of \textit{condition-justified} definitions.

If the previously specified condition is valuable and the kind of equivalence is natural, condition-justification is a reasonable kind of justification.\footnote{For the condition-justified definitions of my case study we will see why the conditions are valuable and the equivalences are natural. Yet
characterising what constitutes valuable conditions or natural kinds of equivalence \textit{at a general level} would require further research.} A condition-justified definition can be regarded as providing knowledge because it answers the question of which definition corresponds naturally to a previously specified condition.

Two definitions of the notions of randomness in ergodic theory [cf.~section~\ref{CS}] are condition-justified, namely weak mixing (all versions) and Bernoulli-system (one version) are condition-justified. Let us discuss them now.

\subsubsection{Weak Mixing}
Recall \textit{Definition~\ref{wmsA}} of weak mixing. Two alternative equivalent definitions are \cite[pp.~65--67]{Petersen1983}:
\begin{definition}\label{wmsII}
$(X,\Sigma,\mu,T)$ is weakly mixing iff for all $A,B\in\Sigma$\begin{displaymath}\lim_{n\rightarrow\infty}\frac{1}{n}\sum_{i=0}^{n-1}|\mu(T^{i}(A)\cap B)-\mu(A)\mu(B)|=0.\end{displaymath}
\end{definition}
\begin{definition}\label{wm}
$(X,\Sigma,\mu,T)$ is weakly mixing iff
for all $f,g\in L^{2}(X,\Sigma,\mu)$\begin{displaymath}\lim_{n\rightarrow \infty}\frac{1}{n}\sum_{i=0}^{n-1}\mid\int f(T^{i}(x))g(x)d\mu-\int f(x)d\mu \int g(x)d\mu\mid=0,\end{displaymath}
\end{definition}
where $L^{2}(X,\Sigma,\mu)$ is the Hilbert space of square integrable functions on $(X,\Sigma,\mu)$.

We already argued that \textit{Definition~\ref{wmsA}} can be natural-world-justified. The first three papers discussing weak mixing seem to be \cite{Hopf1932a}, \cite{Hopf1932b}, and \cite{KoopmanNeumann1932}. These papers show that there is more to say for three reasons.

First, \cite{Hopf1932a} starts by emphasising the importance of ergodicity for statistical mechanics. He then considers another statistical property discussed by Poincar\'{e}: when initially a certain part of a fluid is coloured, experience shows that after a long time the colour uniformly dissolves in the fluid. Mathematically, \citeauthor{Hopf1932a} expresses this by strong mixing.\footnote{$(X,\Sigma,\mu,T)$ is strongly mixing iff for all $A,B\in\Sigma$
\begin{displaymath}
\lim_{n\rightarrow\infty}\mu(T^{n}(A)\cap B)=\mu(A)\mu(B).
\end{displaymath}}
Interested in the interrelationship between strong mixing and ergodicity, he conjectures that a continuous system $(X,\Sigma,\mu,T_{t})$ is strongly mixing iff for all $t$ the discrete-time transformation $T^{t}$ is ergodic.\footnote{A system $(X,\Sigma,\mu,T)$ is \textit{ergodic} iff for all $A\in\Sigma$ if $T(A)=A$, then $\mu(A)=0$ or $1$.} Yet he is unable to prove this (it was later shown to be false). As a result, Hopf attends to the question which weaker statistical property is equivalent to the condition that for all $t$ the transformation $T^{t}$ is ergodic. The answer he arrives at is \textit{Definition~\ref{wm}}. Therefore, weak mixing is condition-justified because its justification stems from it being equivalent in a natural way to a condition regarded as valuable.

Second, \cite{Hopf1932b} is concerned with Gibbs fundamental hypothesis that any initial distribution tends toward statistical equilibrium, and he derives several conditions under which this hypothesis holds true.
Within this context, he becomes interested in how properties of a system $(X,\Sigma,\mu,T)$ relate to the composite system $(X\times X,\Sigma\times\Sigma,\mu\times\mu,T\times T)$ comprising two copies of the single system. Because of the importance of ergodicity, it is natural to ask: which property of the single system is equivalent to the composite system being ergodic? \cite{Hopf1932b} provides the answer, namely \textit{Definition~\ref{wm}} of weak mixing. Hence weak mixing is condition-justified as \citet[p.~1022]{Halmos1949} stresses by referring to \textit{Definition~\ref{wmsA}} and \textit{\ref{wmsII}}: an ``indication that weak mixing is more than an analytic artificiality is in the assertion that $T$ is weakly mixing if and only if its direct product with itself is indecomposable [ergodic]''.

Third, when discussing \textit{Definition~\ref{KSE}}, we encountered the property of a continuous spectrum which arises in spectral theory. \cite{KoopmanNeumann1932} emphasise the naturalness of, and devote their paper to, this property.
From the beginning of ergodic theory the correspondence of concepts from spectral theory and set-theoretic and integral-theoretic concepts from ergodic theory has been a core theme. Hence it was natural to address the question, as \citeauthor{KoopmanNeumann1932} did, which set-theoretic or integral-theoretic definition is equivalent to having a continuous spectrum. The answer they arrived at is \textit{Definition~\ref{wmsA}} of weak mixing. Thus, again, weak mixing is condition-justified.

I have found no book motivating weak mixing by the condition that for all $t$ the transformation $T^{t}$ is ergodic. This might be because that characterisation does not hold for discrete systems.\footnote{The irrational rotation of the circle is a counterexample \cite[p.~8]{Petersen1983}.} The other two interpretations of weak mixing as condition-justified appear in standard books on ergodic theory, e.g.~\citet[p.~39]{Halmos1956} and \citet[p.~64]{Petersen1983}. The latter comments:
\begin{quote}
That the concept of weak mixing is natural and important can be seen from the following theorem, according to which a transformation is weakly mixing if and only if its only measurable eigenfunctions are the constants.
\end{quote}

To summarise, all versions of weak mixing are condition-justified because their justification stems from their being equivalent in a natural way to a condition regarded as valuable. The next definition illustrates the danger of not appreciating that a definition is condition-justified.

\subsubsection{Discrete-time Bernoulli-system}

Recall \textit{Definition~\ref{BerSys}} of a Bernoulli-system. The appeal to isomorphisms makes this definition indirect. Furthermore, most states of the systems encountered in the sciences, e.g.~states of Newtonian systems, are not infinite sequences. Thus it is often easier to work without notions referring to infinite sequences. In investigating simple systems isomorphic to Bernoulli-shifts, it became clear that proving an isomorphism amounts to finding a partition which can be used to code the dynamics. Hence it was natural to ask which condition that does not appeal to isomorphisms and infinite sequences, but to partitions, is equivalent to a Bernoulli-system. \cite{Ornstein1970} gives the definition which answers this question:\footnote{Strictly speaking, it is equivalent to \textit{Definition~\ref{BerSys}} only for Lebesgue spaces. Since all spaces of interest are Lebesgue, this is considered unproblematic \cite[pp.~16--17 and p.~275]{Petersen1983}.}
\begin{definition}\label{BsystemP}
$(X,\Sigma,\mu,T)$ is a Bernoulli system iff there is a partition $\alpha$ such that \\
(i) $T^{i}\alpha$ is
an independent sequence, i.e.\ for any distinct $i_{1},\ldots,i_{r}\in
\field{Z}$, and not necessarily distinct $\alpha_{j}\in\alpha$,
$j=1,\ldots,r$ ($r\geq 1$):\begin{displaymath}
\mu(T^{i_{1}}\alpha_{1}\cap\ldots\cap T^{i_{r}}\alpha_{r})=\mu(\alpha_{1})\ldots\mu(\alpha_{r}).\end{displaymath}
(ii) $\Sigma$ is generated by $\{T^{i}\alpha\,|i\in\field{Z}\}$.
\end{definition}
Hence \textit{Definition~\ref{BsystemP}} can be justified by the fact that it gives an answer to the above question, i.e.\ it is condition-justified. Standard books on ergodic theory also hint at that justification [\citeauthor{Shields1973}, \citeyear{Shields1973}, p.~8, p.~11; \citeauthor{Sinai2000}, \citeyear{Sinai2000}, p.~47].

There have been attempts to justify \textit{Definition~\ref{BsystemP}} as capturing a preformal idea of randomness. Interpreting the measure as time-independent probability, condition (i) captures the idea that any number of finite events of a specific partition at different times are independent. \citeauthor{Berkovitzetal2006}~[\citeyear{Berkovitzetal2006}] argue that because condition (i) can be thus interpreted, Bernoulli-systems capture randomness;\footnote{Actually, a slip occurred in  \citeauthor{Berkovitzetal2006}'s~[\citeyear{Berkovitzetal2006}, p.~667] interpretation of condition (i); (i) holds only for any finite number of events of a \textit{specific} partition at different times, not for \textit{any} events.}
they do not say anything about condition (ii). Yet since (i) is only one part of this definition, this justification of \textit{Definition~\ref{BsystemP}} fails.\footnote{For instance, the following system fulfills (i) but not (ii):  let $(X,\Sigma,\mu)$ be the ordinary Lebesgue measure space of the unit cube $X$. Let\begin{displaymath}
T(x,y,z):=(2x,\frac{y}{2},z)\,\,\textnormal{if}\,\,0\leq x<\frac{1}{2},\,\,(2x-1,\frac{y+1}{2},z)\,\,\textnormal{if}\,\,\frac{1}{2}\leq x\leq 1.\end{displaymath}
Obviously, for $(X,\Sigma,\mu,T)$ condition (i) of \textit{Definition~\ref{BsystemP}} holds for $\alpha=\{\{x\in X\,|$\linebreak $\,0\leq x<\frac{1}{2}\},\{x\in X\,|\,\frac{1}{2}\leq x\leq 1\}\}.$ But $(X,\Sigma,\mu,T)$ is not a Bernoulli-system since it is not ergodic.}   Generally, if a definition does not capture the idea it is said to capture, the justification fails because it is unclear why this definition is chosen.

\citeauthor{Batterman1991}'s~[1991] and \citeauthor{Sklar1993}'s~[\citeyear{Sklar1993}, pp.~238--239] motivation for \textit{Definition~\ref{BsystemP}} is also that it captures a preformal idea of randomness. Their argument as expressed by \citeauthor{Batterman1991}~[\citeyear{Batterman1991}, pp.~249--250] is:
\begin{quote}
Now let us see just how random a Bernoulli system is [...] The Bernoulli systems are those in which knowing the entire past history of box-occupations even relative to a partition that is generating in the above sense, is insufficient (in the sense of being probabilistically independent) for improving the odds that the system will next be found in a given box.
\end{quote}
As an interpretation of randomness this is puzzling. Even if it exactly corresponded to \textit{Definition~\ref{BsystemP}},\footnote{It does not. First, their interpretation does not make clear that the matter of concern is the \textit{existence} of a partition satisfying (i) and (ii). Even if this is disregarded, their interpretation applies to more systems than Bernoulli systems, for instance, also for to the system and the partition~$\alpha$ discussed in the previous footnote. Here the events constituting any entire history of box-occupations are of probability zero. Then \citeauthor{Batterman1991}'s and \citeauthor{Sklar1993}'s claim about the independence from the entire past history of box-occupations is trivially true. Correct is: any finite number of events of a specific partition at different times are independent, even though the partition is generating.}
it is unclear, \textit{from the viewpoint of capturing a preformal idea of randomness}, why independence is required relative to \textit{generating} partitions; and I found no convincing justification for this.

It seems that the difficulty stems from the fact that \textit{Definition~\ref{BsystemP}} is really condition-justified. As we have seen for weak mixing, condition-justified definitions may in other contexts also capture a preformal idea valuable in some sense. However, often---and this is true for \textit{Definition~\ref{BsystemP}} as discussed---this won't be the case. Then there is the danger of not appreciating that a definition is condition-justified and claiming that it captures a valuable preformal idea, when it does not. It seems that in interpreting \textit{Definition~\ref{BsystemP}} \citeauthor{Batterman1991} and \citeauthor{Sklar1993} fell into this trap. This danger is similar to the one identified by \cite[p.~153]{Lakatos1976}, viz.\ claiming that a proof-generated definition captures a valuable preformal idea when it does not.

Let us now turn to the final kind of justification I have identified.

\subsection{Redundancy-Justification}
We call a \textit{definition which is justified because it eliminates at least one redundant condition in an already accepted definition} \textit{redundancy-justified}. A redundancy-justified definition can be regarded as providing knowledge since it shows that specific conditions in an accepted definition are redundant.

It is obviously desirable in mathematics to find out whether there are any redundant conditions in an already accepted definition. Typically, both the original definition, and the one in which the redundant conditions are eliminated, each have their own advantages. It depends on the definitions, but the former might be easier to understand or might allow for a more fine-grained analysis; the latter is simpler (in the sense of being more concise), and it might be that only the latter is easier to use in proofs, allows for natural generalizations, or suggests important analogies.

So when is it better to propound the original definition? And when is it better to introduce instead the new definition without the redundant conditions, i.e.\ when is redundancy-justification a reasonable kind of justification? I think the answer depends on the definition and the context in which the definition is considered.
For the purpose of an introductory textbook it might be better to propound the original definition because it is easier to understand. Conversely, for the purpose of a research article it might be better instead to use the new, concise definition, since it is easier to use in some proofs. Furthermore, in many cases it does not seem to matter much whether the original definition or the definition in which the redundant conditions are eliminated is introduced, so long as the origin of the definition and the redundant conditions are clearly pointed out.

As in the case of proof-generation and condition-justification there is the danger of  not understanding that a definition is redundancy-justified and claiming that it captures a valuable preformal idea, when it does not.

Two definitions of the list of notions of randomness [cf.~section \ref{CS}] are redundancy-justified: the continuous version of a Bernoulli-system, which we will discuss for illustration, and a Kolmogorov-system [\citeauthor{Sinai1963}, \citeyear{Sinai1963}, pp.~64--65; \citeauthor{Uffink2006}, \citeyear{Uffink2006}, pp.~94--96].

\subsubsection{Continuous-time Bernoulli-system}
We have seen that \cite{Kolmogorov1958} and \cite{Sinai1959} established that isomorphic discrete-time Bernoulli-systems have the same Kolmogorov-Sinai entropy [cf.~subsection~\ref{KSentropy}]. A decade later \cite{Ornstein1970} proved the converse, i.e.~that Bernoulli-systems with equal entropy are isomorphic.

Having established that celebrated result, Ornstein became interested in finding an analogous definition of a Bernoulli-system for continuous time, and he asked whether the Kolmogorov-Sinai entropy could be used to classify them too. The most obvious definition of a continuous system $(X,\Sigma,\mu,T_{t})$ describing an independent process is that for all $t$ the discrete system $(X,\Sigma,\mu, T_{t})$ is a Bernoulli-system. \cite{Ornstein1973} first introduces this definition of a continuous Bernoulli-system, and then he shows that there are redundant conditions in this definition because it is equivalent to the following definition:
\begin{definition}\label{cb}
$(X,\Sigma,\mu,T_{t})$ is a Bernoulli-system iff the discrete system  $(X,\Sigma,\mu,T_{1})$ is a Bernoulli-system.
\end{definition}
Hence \textit{Definition~\ref{cb}} is redundancy-justified because it eliminates redundant conditions. In this way it seems to be justified in Ornstein's [\citeyear[p.~56]{Ornstein1974}] book too.\footnote{Ornstein [\citeyear[p.~56]{Ornstein1974}] expresses this indirectly by introducing continuous Bernoulli systems as follows; `We will call a flow \{$(X,\Sigma,\mu,T_{t})$\} a \{continuous Bernoulli-system\} if \{$(X,\Sigma,\mu,T_{1})$\} is a \{Bernoulli system\}. (We will prove later that if \{$(X,\Sigma,\mu,T_{1})$\} is a \{continuous Bernoulli system\}, then \{$(X,\Sigma,\mu,T_{t_{0}})]$\} for each fixed $t_{0}$ is a \{Bernoulli system\}').}

Any continuous Bernoulli-system can be normalised such that it has Kolmogorov-Sinai entropy one. \cite{Ornstein1973} indeed showed that any (and only a) normalised continuous-time Bernoulli-system is isomorphic to a normalised continuous-time Bernoulli-system. Later it was proven that some Newtonian systems are Bernoulli, showing how random Newtonian systems can be \cite[p.~183]{Ornstein1989}.

\subsection{Occurrence of the Kinds of Justification}\label{subject}
To sum up: in addition to Lakatos's proof-generated definitions, I have identified three kinds of justification of definitions. To my knowledge, condition-justification and redundancy-justification have not been identified before. I do not claim that the kinds of justification we have discussed are the only ones at work in mathematics. Further studies might unveil yet other ones.

Two more general comments about justifying definitions should be added here. First, for any kind of justification there are three possibilities: (i) a definition is reasonably justified in this way; (ii) it is justified but not reasonably justified in this way; (iii) it is not justified in this way. As regards (ii), for instance, if the idea of being equivalent in a measure-theoretic sense to an independent process like throwing a die was not valuable, \textit{Definition~\ref{BerSys}} would be natural-world-justified but not reasonably justified. Second, an already justified definition  has sometimes additional good features which support this definition but which do not by themselves constitute a sufficient justification. These features may also be important in deciding between different definitions. For instance, it is often said that  a merit of the Kolmogorov-Sinai entropy is its neat connection to other notions of randomness like Kolmogorov-systems. These are good features but not sufficient justifications; since if there were no further reasons for studying the definition, there would still be the question why we should regard it as worth considering [cf.~\citeauthor{Smith1998}, \citeyear{Smith1998}, pp.~174--175].

How widely do the kinds of justification we discussed occur? To answer this, I first comment on the notion of a mathematical subject. I think that regardless of which plausible understanding of `subject' is adopted, my claims are true. But a possible way to operationalise this idea is the following: with the subjects identified by the  \textit{Mathematical Subject Classification}\footnote{This is a five digits
classification scheme of subjects formulated by the American Mathematical Society; see~www.ams.org/msc. For our purposes subjects concerned with education, history or experimental studies have to be excluded.} it would be possible to create a list of subjects of the mathematics from the nineteenth  century up to today. Then with notions of randomness in ergodic theory some of the main definitions of the subject `strange attractors, chaotic dynamics' have been investigated.

Based on my knowledge of mathematics, I endorse the following claims about mathematics produced in the twentieth century and up to the present day:\footnote{Starting with the twentieth century is somewhat arbitrary. All the here-discussed kinds of justification appear also important in nineteenth century mathematics. Yet older mathematics may be significantly different. Hence a close investigation would be necessary to identify the role the kinds justification play in older mathematics.} \textit{all the kinds of justifications I have discussed in this paper are widespread.} More specifically, proof-generated, condition-justified, and redundancy-justified definitions are all found in the majority of mathematical subjects with explicit definitions. Also, for nearly all mathematical subjects with explicit definitions which (among other things) aim at describing or understanding the natural world, natural-world-justified definitions are found. This includes subjects not only from what is called applied mathematics but also from pure mathematics, e.g., measure theory. Furthermore, as in my case study, \textit{for nearly all mathematical subjects with explicit definitions many different ways of justifying definitions are found and are reasonable.} Indeed, I would be surprised if one subject could be found where only one kind of justification is important. Clearly, my case study shows that for the subject `strange attractors, chaotic dynamics' these claims hold true.

For my case study the argumentation involved in justifying definitions is typically not explicitly stated but is merely hinted at or merely implicit in the mathematics. Because of the conventional style of mathematical writing, this appears to be generally the case in mathematics, as also \citeauthor{Lakatos1976}~[\citeyear{Lakatos1976}, pp.~142--144] claimed.
As we have seen, awareness of the ways of justifying definitions is important for understanding mathematics and for preventing mistakes. Thus it would be desirable if publications addressed this issue more explicitly. Also, it should be mentioned that detailed knowledge of parts of ergodic theory is necessary to assess how definitions are justified in my case study. This confirms \citeauthor{Tappenden2008a}'s claim that judgments about definitions require detailed knowledge of the relevant mathematics [cf.~section~\ref{SOTA}].

Let us reflect on the interrelationships between the kinds of justification, an issue which seems not discussed in the literature.

\section{Interrelationships Between the Kinds of Justification}\label{interr}

In what follows when we speak of an argument for a definition we mean that a reason is provided for a definition which cannot be split into two separate reasons for this definition. Now we first ask about the \textit{interrelationships in one argument}: assume that a specific argument establishes that a definition is justified according to one kind of justification. Can it be that this argument implies that the definition is at the same time also justified according to another kind of justification? Intuitively, one might think that in an argument a definition can only be justified according to one kind of justification. Yet, as we will see, the matter is more complicated. Second, we ask about the \textit{interrelationships between the kinds of justification in different arguments}: if different arguments justify the same definition, what combination of kinds of justification do we find? We will discuss these two cases in the next two subsections.

\subsection{One Argument}

Clearly, there are arguments where a definition is only proof-justified, natural-world-justified, condition-justified or redundancy-justified. For example, uniform convergence as discussed by \citeauthor{Lakatos1976}~[\citeyear{Lakatos1976}, pp.~131--133] is only proof-justified, \textit{Definition~\ref{BerSys}} of a Bernoulli-system as capturing the idea of a measure-theoretic system being equivalent to an independent process is only natural-world-justified, weak mixing as corresponding to ergodicity of the composite system is only condition-justified, and \textit{Definition~\ref{cb}} of a continuous-time Bernoulli system as eliminating redundant conditions is only redundancy-justified.

By going back to the characterisation of the kinds of justification, we see that the intuition that in an argument a definition can only be (reasonably) justified according to one kind of justification is correct except for one case. Namely, in rare cases \textit{condition-justified definitions are at the same time proof-generated} in an argument. This is so if and only if the kind of equivalence is regarded as natural because it occurs in the formulation of a conjecture that should be established. For example, assume the following conjecture is regarded as valuable: each function in a convergent sequence of functions is continuous if and only if the limit function of the convergent sequence is continuous. Further, assume that sequences of pointwise convergent continuous functions without continuous limit functions are known. Then mathematicians might ask: how has the notion of convergence to be changed such that if and only if the limit function is continuous the sequence of continuous functions is convergent? The definition answering this question would be clearly condition-justified. But it would also be proof-generated since it is needed in order to prove the above conjecture.

Let us now turn to the interrelationships in different arguments.

\subsection{Different Arguments}
In our case study different arguments establish that weak mixing is condition-justified: weak mixing corresponds to ergodicity of the composite system, to the set-theoretic or integral-theoretic condition equivalent to having a continuous spectrum, and for continuous systems to the condition that for all $t$ the transformation $T^{t}$ is ergodic. Generally, \textit{one and the same definition can be (reasonably) justified in the same way in different arguments by referring to different conjectures, preformal ideas etc}. For proof-generation \citeauthor{Lakatos1976}~[\citeyear{Lakatos1976}, pp.~127--128] recognises this pattern.

What is more, we have seen that in different arguments \textit{Definition~\ref{wmsA}} of weak mixing is justified in different ways:
as mentioned above, it is condition-justified but also natural-world-justified, expressing the idea that almost all sufficiently past events are approximately independent. Likewise, the Kolmogorov-Sinai entropy is natural-world-justified, expressing the idea of the highest average amount of information produced per step relative to a coding; but it is also proof-generated concerning the conjecture that spectrally-isomorphic systems are not always isomorphic. Generally, \textit{one and the same definition can in different arguments be (reasonably) justified in different ways}.

Finally, a definition which is justified in any way can be used to (reasonably) justify a definition in an arbitrary way. In this sense the different kinds of justification are closely connected. For example, the natural-world-justified \textit{Definition~\ref{BerSys}} of a Bernoulli-system is used to justify the condition-justified \textit{Definition~\ref{BsystemP}} of a Bernoulli-system.

A special case of this is when for proof-generated definitions preformal-ideas shine through (which can be, but does not have to be the case). For instance, consider definitions of polyhedron as discussed by \citeauthor{Lakatos1976} [\citeyear{Lakatos1976}]. Early definitions of polyhedron, which seem to be justified because they capture the preformal idea of a solid with plane faces and straight edges, were eventually replaced by definitions which are needed to prove the Euler conjecture. For these proof-generated definitions, to some extent, the preformal idea of the old definitions still shine through. Hence \citeauthor{Lakatos1976}'s~[\citeyear{Lakatos1976}, p.~90] claim ``In the different proof-generated theorems we have nothing of the naive concept''  is an unfortunate exaggeration.

We now return to Lakatos's ideas on justifying  definitions.

\section{Assessment of Lakatos's Ideas on Proof-Generated Definitions}\label{Lakatos}

First, in focusing on proof-generated definitions, \textit{Lakatos fails to recognise the interplay between the different kinds of justification of definitions}, which we discussed in section~\ref{interr}. In particular, Lakatos never indicates that in different arguments the same definition can be justified in different ways.

Second, \textit{Lakatos fails to show, as we did for notions of randomness in ergodic theory, that often various kinds of justification can be found and that a variety of kinds of justification can be reasonable}.  We argued that Lakatos may have believed the following [cf.~section \ref{SOTA}]: there are many mathematical subjects where proof-generation should be the sole important way that definitions are justified; and there are many subjects after the discovery of proof-generation where proof-generation is the sole important way that definitions are justified. From our claim that for nearly all mathematical subjects many different ways of justifying definitions are found and are reasonable follows that this must be wrong [cf.~section \ref{subject}]. That is, subjects created after the discovery of proof-generation where solely proof-generated definitions are found and are reasonable appear to be exceptional.

Indeed, Lakatos could have shown with his case studies that often various kinds of justification are found and that various kinds of justification can be reasonable. To demonstrate this, I will now show that even for the subjects discussed by \citeauthor{Lakatos1976}~[\citeyear{Lakatos1976}] not only proof-generation but also other kinds of justification are important. Because of lack of space, I show this here only for the subjects to which the definition of uniform convergence and the \citeauthor{Caratheodory1914} definition of measurable sets belong. But this hypothesis can easily seen to be also true for the subjects to which the other proof-justified definitions Lakatos discusses (namely polyhedron, bounded variation and the Riemann integral) belong.

\citeauthor{Lakatos1976}~[\citeyear{Lakatos1976}, pp.~144--146] argues that \textit{uniform convergence} is proof-generated, also by referring to textbooks.  This definition falls under the subject `convergence and divergence of series and sequences of functions'.\footnote{In the terms of the Mathematical Subject Classification.} A definition discussed in this subject is the radius of convergence of a power series. A power series is of the form $\sum_{k=0}^{\infty}a_{k}(x-x_{0})^{k}$, where $a_{k}\in\field{R}$.
\begin{definition}\label{ende2}
Its radius of convergence is the unique number $R\in[0,\infty]$ such that the series converges absolutely if $|x-x_{0}|<R$ and diverges if $|x-x_{0}|>R$.
\end{definition}

The radius of convergence is often defined differently as follows. The \textit{root test} is a powerful criterion for the convergence of infinite series. Hence the question arises whether there is a definition which is equivalent to the radius of convergence as defined above but which gives an explicit way to calculate this radius by referring to the root test. The answer is yes, namely:
\begin{definition}\label{ende}
For a power series the radius of convergence is\begin{displaymath}
R:=1/\limsup_{k\rightarrow \infty}\sqrt[k]{|a_{k}|}.\end{displaymath}
\end{definition}
Thus \textit{Definition~\ref{ende}} is condition-justified, as, for example, hinted at in Marsden and Hoffman's [\citeyear{MarsdenHoffman1974}, pp.~289--290] standard analysis book:
``The reason for the terminology in \{\textit{Definition~\ref{ende}}\} is brought out by the following result [that by applying the root test, \textit{Definition~\ref{ende}} is equivalent to \textit{Definition~\ref{ende2}}].''

\citeauthor{Lakatos1976}~[\citeyear{Lakatos1976}, pp.~152--154], mainly by referring to Halmos's [\citeyear{Halmos1950}] book, argues that the \textit{\citeauthor{Caratheodory1914} definition of measurable sets} is proof-generated. This definition falls under the subject `classes of sets in measure theory.\footnote{`Classes of sets, measurable sets, Suslin sets, analytic sets' in the terms of the Mathematical Subject Classification.} The definition of a $\sigma$-algebra clearly belongs to this subject. The basic idea of a $\sigma$-algebra is to have a collection of subsets of $X$ including $X$ which is closed under countable set-theoretic operations. Thus a usual definition is \cite[pp.~1--2]{Cohn1980}:
\begin{definition}
A set $\Sigma$ of subsets of $X$ is a $\sigma$-algebra iff\\
(i) $X\in\Sigma$,\\
(ii) for all $A\subseteq X$ if $A\in\Sigma,\,\,\textnormal{then}\,\,X\setminus A\in\Sigma$,\\
(iii) for all sequences $(A_{k})_{k\geq 0}$ if $A_{k}\in\Sigma$ for all $k\geq 0,\,\textnormal{then}\,\bigcup_{i=0}^{\infty}A_{k}\!\in\!\Sigma$,\\
(iv) for all sequences $(A_{k})_{k\geq 0}$ if $A_{k}\in\Sigma$ for all $k\geq 0,\,\textnormal{then}\,\bigcap_{i=0}^{\infty}A_{k}\!\in\!\Sigma$.
\end{definition}
Now one can easily see that the conditions (i), (ii) and (iii) imply (iv). Consequently, many use the following definition because it eliminates a redundant condition.
\begin{definition}
A set $\Sigma$ of subsets of a set $X$ is a $\sigma$-algebra iff (i), (ii) and (iii) hold.
\end{definition}
Clearly, it is redundancy-justified as, for instance, in Ash's [\citeyear[p.~4]{Ash1972}] standard book on measure theory.

To conclude, even for the subjects discussed by Lakatos various kinds of justification are found and are reasonable.

\section{Conclusion}
This paper addressed the actual practice of how definitions in mathematics are justified in articles and books and whether the justification is reasonable. In section \ref{SOTA} I discussed the main account of these issues, namely Lakatos's ideas on proof-generated definitions. While important, this paper showed how they are limited. My assessment of Lakatos and my thoughts on justifying definitions are based on a case study of notions of randomness in ergodic theory, which was introduced in section \ref{CS}.
In section \ref{chaos} I identified three other important and common ways of justifying definitions: natural-world-justification, condition-justification and redundancy-justification. To my knowledge, condition-justification and redundancy-justification have not been discussed so far. Also, we have seen that awareness of the ways definitions are justified is important for mathematical understanding and for avoiding mistakes.
Then in section \ref{interr} I discussed the interrelationships between the different kinds of justification of definitions, an issue which has not been addressed before. In particular, I argued that in different arguments the same definition can be justified in different ways. Finally, in section \ref{Lakatos} I criticised Lakatos's ideas on proof-generated definitions. They fails to recognise the interplay between the kinds of justification. Furthermore, they fail to show that often various kinds of justification are found and that a variety of kinds of justification can be reasonable.

\section*{Acknowledgments.}
I am indebted to Jeremy Butterfield and Peter Smith for valuable feedback on previous versions of this manuscript. Many thanks to Roman Frigg, Franz Huber, Brendan Larvor, Mary Leng, Paul Weingartner, two anonymous referees and the audiences at the 1st London-Paris-Tilburg Workshop in Logic and Philosophy of Science and the 1st Conference of the European Philosophy of Science Association for helpful comments. I am grateful to St John's College, Cambridge, for financial support.

\end{document}